
\input amstex
\documentstyle{amsppt}
 \magnification=1200 
\voffset=-0.6in
\vsize=7.5in

\def\e{\epsilon}
\def\D{\Cal D}
\def\ol{\overline}
\def\lf{\left}
\def\ri{\right}

\def\BY{\frak m _{\text{BY}}}

\def\wt{\widetilde}

\def\p{\partial}

\def\R{\Bbb R}
\def\cR{\Cal R}
\def\lgle{\langle}
\def\rgle{\rangle}

\document

\topmatter

\leftheadtext{ Yuguang Shi and Luen-fai Tam} \rightheadtext{Some
lower estimates of ADM mass and Brown-York mass}
  \document
\topmatter
\title{Some lower estimates of ADM mass and Brown-York mass}\endtitle
\author{Yuguang Shi\footnotemark and Luen-Fai Tam\footnotemark}\endauthor
 \footnotetext"$^1$"{Research partially supported by NSF of China, Project number 10001001.}
\footnotetext"$^2$"{Research partially supported by   Earmarked
Grant of Hong Kong \#CUHK4032/02P.}
\address{Key Laboratory of  Pure and Applied Mathematics, School of Mathematics Science,
Peking University, Beijing, 100871, China}
\endaddress
\email{ygshi\@math.pku.edu.cn}
\endemail
\address
{Department of Mathematics, The Chinese University of Hong Kong,
Shatin, Hong Kong, China}
\endaddress
\email{lftam\@math.cuhk.edu.hk}
\endemail
\affil {
Peking University\\
 The Chinese University of Hong Kong,
 }
\endaffil
\date June, 2004\enddate
\abstract
We give some lower estimates of the ADM mass of an asymptotically flat (AF) Riemannian manifold without assuming that the scalar curvature of the manifold is nonnegative. Some  sufficient conditions for an AF   manifold to have nonnegative ADM mass are obtained. We  also give some lower estimates of the Brown-York mass of a compact three manifold with smooth boundary. From these estimates, we   generalize some previous results of the authors.
\endabstract
\endtopmatter

\subheading{\S1 Introduction}\vskip .2cm

More than twenty years ago, Schoen and Yau \cite{SY1-2}  proved
the positive mass theorem. Later, using spinors Witten \cite{W}
gave a simple proof of the result. A mathematical rigorous proof
of Witten's argument was  given by Parker and Taubes \cite{PT},
see also \cite{B1}. For the time-symmetric case, the positive mass
theorem asserts that the Arnowitt-Deser-Misner (ADM) mass of each
end of a three dimensional asymptotically flat (AF) manifold $M$
with finitely many ends with $L^1$ integrable and nonnegative
scalar curvature  is nonnegative. Moreover,  if  the ADM  mass of one of
the ends is zero  then   the manifold has only one end and is isometric to the
three dimensional Euclidean space. In the time-symmetric case, the scalar curvature is nonnegative
means physically  that the local mass density is nonnegative. The
condition that the scalar curvature being in $L^1$ is necessary in
order that the ADM mass is well-defined, see \cite{B1}.

In \cite{ZZ}, L. Zhang and X. Zhang studied an interesting
question asked by S.-T. Yau how the condition that the scalar
curvature is nonnegative can be relaxed so that the ADM mass of an
AF manifold is still nonnegative. In \cite{ZZ}, they proved that
the positive mass theorem is still true under the assumptions that
the first Dirichlet  eigenvalue and the first eigenvalue of
Neumann type of the conformal Laplacian operator are nonnegative.
Motivated by the above mentioned results, in this work we shall
discuss lower bounds of the ADM mass of an AF manifold without
assuming that the scalar curvature is nonnegative. See Theorems
3.1 and 3.2 for more details. From these results, we obtain
conditions in terms of the geometry of the AF manifold so that the
ADM mass of the manifold is nonnegative. More precisely, we have
the following (Corollaries 3.1 and 3.3):

\proclaim{Theorem 1.1} Let $(M^3,g)$ be an AF manifold with one
end such that its scalar curvature $\cR$ is in $L^1(M)$. Let
$\cR_+$ and $\cR_-$ be  positive part and negative part of  $\cR$
respectively. Suppose one of the following   is true, then the ADM
mass of $M$ is nonnegative. \roster

\item"{(a)}"
$
\lf(\int_M\lf(\frac{\cR_-}4\ri)^{\frac{3}{2}}\ri)^{\frac{2}{3}}
<\frac\Lambda2, $
and
$$
\lf[\Lambda-2\lf(\int_M\lf(\frac{\cR_-}4\ri)^{\frac{3}{2}}\ri)^{\frac{2}{3}}\ri]\int_M\cR_+\ge
\lf[\Lambda
+2\lf(\int_M\lf(\frac{\cR_+}4\ri)^{\frac{3}{2}}\ri)^{\frac{2}{3}}\ri]\int_M\cR_-.
$$
where $\Lambda$ is the Sobolev constant of $(M,g)$.

\item"{(b)}" $
\lf(\int_M\lf(\frac{\cR_-}4\ri)^{\frac{3}{2}}\ri)^{\frac{2}{3}} <
\Lambda, $ and

$$
B\ge C_1A^2\lf[\sup_M|Rm|\int_MR_-+||\nabla Rm||_2\lf(
\int_MR_-\ri)^\frac12  \ri]
$$
where 
$$
B=\inf_{D}\lf\{\int_{M\setminus D}|Rm|^2\bigg|\ \ \text{vol}(D)\le
C_2\lf[ A\Lambda^{-1} \int_M\cR_-  \ri]^{3}\ri\},
$$
$C_1$ and  $C_2$ are   positive absolute constants, and $A$ is a positive constant depending only on $\Lambda$ and $\cR_-$.
\endroster
\endproclaim

A condition similar to (a) can also be obtained from the proof of
\cite{ZZ}. Condition (a) implies that $\int_M\cR\ge0$. In fact the
assumptions  in \cite{ZZ} that the first eigenvalue of Neumann
type is nonnegative also imply that $\int_M\cR\ge0$. Hence if
$\int_M\cR<0$, then we cannot apply (a) of the theorem.

In \cite{BF, FK}, Bray-Finster and Finster-Kath  obtained  lower
estimates of the ADM mass in terms of the curvature under the
assumption that the scalar curvature is nonnegative. Condition (b)
in Theorem 1.1 is obtained by generalizing their results to the
case that the scalar curvature may be negative somewhere. If $\cR$
is nonnegative, then $B\geq 0$, condition (b) is automatically
satisfied.

In the second part of this work, we shall discuss similar issues
for the Brown-York mass of a compact manifold with smooth
boundary. 

Let  $(\Omega^3,g)$ be a $3$-dimensional compact manifold with
smooth boundary. For simplicity, in this paper we always assume
that $\p\Omega=\Sigma$ is connected. Suppose the   Gauss curvature
of $\Sigma$ is positive, then by a classical result we know that
$\Sigma$ can be isometrically embedded into $\Bbb{R}^3$.  Let
$H_0$ be the mean curvature  of image of the embedding in
$\Bbb{R}^3$ with respect to the  outward  norm.  The Brown-York
mass of $(\Omega,g)$ is defined as   (see \cite{BY 1-2}):

$$\BY(\Omega)= \int_\Sigma (H_0 -H)d\sigma$$
where $d\sigma$ is the volume element of $\Sigma$, and $H$ is the mean
curvature of $\p \Omega$ with respect to original metric $g$ and
outward   norm. In our convention, the mean curvature of the unit
sphere in $\R^3$ is 2. In \cite{ST1}, under the assumption that
$\Sigma$ has positive Gauss curvature, it was proved
 that if (i) the scalar curvature of $(\Omega,g)$
is nonnegative and (ii) $H>0$,  then $\BY(\Omega) \geq 0$.
Moreover, equality holds if and only if $\Omega$ is a domain in
$\R^3$. The condition on the Gauss curvature of $\Sigma$   has
been relaxed in \cite{ST2} where the Gauss curvature of the
boundary is only assumed to be nonnegative. In the second part of
this paper, we shall try to relax condition (i) or condition (ii).
As for ADM mass of an AF manifold, we shall first    give some
lower estimates of the Brown-York mass.

Consider a compact manifold $(\Omega^3,g)$ with smooth boundary
$\Sigma$ which has positive Gauss curvature. Let $\cR$ be the
scalar curvature of $M$   and let $s_0>0$ be such that
$d(x,\p\Omega)$ is smooth in $\Omega_{s_0}=\{x|\
0<d(x,\p\Omega)\le s_0\}$. We have the following (see Theorem  4.1
and Corollary  4.2):
\proclaim{Theorem 1.2} With the above notations, suppose the mean
curvature of the level set $\{x|\ d(x,\Sigma)=s\}$ is positive
with respect to the outward normal for $0\le s\le s_0$. Then there
exists a constant $C>0$ depending only on $s_0$, $H_{\min}$,
$\Lambda$ and $|\Omega|$ where $H$ is the mean curvature of
$\Sigma$ and $\Lambda$ is the Sobolev constant of $\Omega$ such
that if

$$
\sup_\Omega \cR_-\le C,
$$
then the Brown-York mass of $\Omega$ satisfies:
\roster

\item"{(a)}"
$$
\BY(\Omega)\ge\frac{\Lambda -2\beta}{32(\Lambda
+\delta-\beta)}\lf(\int_\Omega
\cR_+-\frac{\Lambda+2\delta}{\Lambda-2\beta}\int_\Omega\cR_-\ri)\tag1.1
$$
provided that $ \beta=\lf(\int_\Omega \frac18\cR_-^{\frac
32}\ri)^\frac23< \frac{\Lambda}2 $
where $\delta=\lf(\int_\Omega \frac18\cR_+^{\frac 32}\ri)^\frac23$; and

\item"{(b)}"
$$
\BY(\Omega)\ge\frac{\lambda}4\int_\Omega
\frac{\cR}{8\lambda+\cR}\tag1.2
$$
provided that  $8\lambda+\cR>0$ in $\Omega$ where $\lambda$ is the
first Dirichlet eigenvalue of the Laplacian of $\Omega$.
\endroster
\endproclaim

From the theorem one can conclude that the Brown-York mass is still nonnegative if the mean curvature is positive
 and the scalar curvature is not very negative. One may ask what would happen if the mean curvature is negative
somewhere. In this case we have the following, see Theorem 4.2 and
Corollary 4.3: \proclaim{Theorem 1.3} With the same notations as
in Theorem 1.2, suppose $\cR\ge0$ and let $\xi=\frac14
\cR_{\min}^\frac12$   where $\cR_{\min}=\inf_{\Omega_{s_0}}\cR$.
Let  $H$ be  the mean curvature of $\{x|\ d(x, \p\Omega)=s\}$,
$0\le s\le s_0$, $H_+=\max\{H,0\}$ and
$H_{\min}=\min_{\p\Sigma}H$. Suppose

\roster
\item"{(i)}" $\xi\ge H_+\tanh(\xi s_0)$ in $\Omega_{s_0}$;
and
\item"{(ii)}" $\xi\tanh(\xi s_0)\ge -4H_{\min}$.
\endroster

Then
$$
\BY(\Omega)\ge \frac14|\Sigma|\, \xi\tanh(\xi s_0).
$$
Moreover, the $\BY(\Omega)$ can also be estimated from below as in (1.1) and (1.2).
\endproclaim

From Theorem 1.3, one can conclude that the Brown-York mass is still nonnegative  if the scalar curvature
is positive but the mean curvature is not too negative. The conditions of the theorem are satisfied for some $s_0$
if $H\ge0$. In particular, (1.1) and (1.2) give lower bounds for the Brown-York mass for compact manifolds with
nonnegative scalar curvature such that its boundary has positive Gauss curvature and positive mean curvature.

It seems likely that the Brown-York mass is nonnegative
irrespectively to the sign of the mean curvature. For example, if
the mean curvature is negative everywhere, then obviously the
Brown-York mass is nonnegative provided the boundary has positive Gauss curvature. However, it is still unclear whether this is true in general.

The paper is organized as follows. In \S2, we shall give
conditions on the existence of  asymptotically constant harmonic
spinors on AF manifolds. The results will be used in the next section. In \S3, we shall obtain different types of lower bounds of the ADM mass in terms of the geometry of the AF manifold. In \S4, we  shall relax the assumptions on the scalar curvature of a
compact manifold or the mean curvature of its boundary and generalize some results in \cite{ST1}. To do this we shall first obtain lower bounds of the Brown-York mass of a compact three manifold with smooth boundary. In \S5, we shall construct examples
which are related to the results in  Theorem 1.1.  We shall   give some interesting applications of quasi-spherical metrics on the
relation between the  classical Minkowski inequalities for convex
surfaces and the positive mass theorem of Herzlich \cite{H}.

\subheading{\S2 Construction of harmonic spinors}\vskip .2cm

In this section,  we shall discuss conditions on an AF spin manifold (see the definition below)
so that one can construct harmonic spinors which are asymptotically parallel near infinity.
Here we do not assume that the scalar curvature is nonnegative. The construction will be used to
give an expression for the ADM mass of the manifold. With the help of these, we shall estimate the lower bound
of the ADM mass and obtain some conditions so that the ADM mass is nonnegative. For simplicity, we always assume
the AF manifold has only one end.

A   complete noncompact spin manifold $(M^n,g)$, $n\ge 3$, is said to be {\it asymptotically flat (AF)} if
 there is a compact set $K$ and a diffeomorphism $\phi:\R^n\setminus{B_R(0)}\to M\setminus K$ for some Euclidean ball
  $B_R(0)$ with center at the origin, such that in the standard coordinates of $\R^n$,

$$
g_{ij}=\delta_{ij}+b_{ij}
$$
with

$$
||b_{ij}||+r||\p b_{ij}||+r^2||\p\p
b_{ij}||=O\lf(r^{2-n}\ri)\tag2.1
$$
where $r$ and $\p$ denote Euclidean distance and the standard
gradient operator on $\R^n$ and the norms are taken with respect
to the Euclidean metric. Moreover, the scalar curvature $\cR$ of
$M$ is assumed to be in$L^1(M)$ so that  the ADM mass of $M$ is
well-defined by  \cite{B1}. Here the ADM mass of $M$ is given by

$$
c(n)\frak
m=\lim_{r\to\infty}\int_{S(r)}\lf(g_{ii,j}-g_{jj,i}\ri)\,dS^i
$$
where $S(r)$ is the Euclidean sphere of radius $r$, $dS^i$ is the normal surface area of $S(r)$ and $c(n)>0$
is a normalizing constant.

First, we need some results on the existence of positive solutions of equation of the form

$$
Lu=\Delta u-qu=0\tag2.2
$$
where $q=O(r^{-n})$. Here and below, when we say a function $f=O(r^\alpha)$ we mean that $|f(x)|\le C(1+r(x))^\alpha$ for
some constant $C$ for all $x\in M$ and $r(x)$ is the geodesic distance of $x$ from a fixed point.

\proclaim{Lemma 2.1} There is a constant $C$ depending only on $M$ and $q$ such that if  $u$ is a positive
solution of (2.2), then for any $r>0$,

$$
\sup_{B_{2r}(p)\setminus B_p(r)}u\le C\inf_{B_{2r}(p)\setminus
B_p(r)}u
$$
for all $r>0$, where $p$ is a fixed point and $B_p(r)$ is the geodesic ball of radius $r$ with center at $p$.
\endproclaim

\demo{Proof} Since $q=O(r^{-n})$ and $g_{ij}$ is uniformly equivalent to the Euclidean metric, the result
follows from  \cite{GT, Theorem 8.20}, see also  \cite{ZZ, p. 666}.
\enddemo

Suppose the first Dirichlet eigenvalue of $L$ in (2.2) is nonnegative, namely:

$$
\int_M\lf(|\nabla v|^2+qv^2\ri)\ge0
$$
for all $v\in C_0^\infty(M)$. By a well-known fact, 
the first Dirichlet eigenvalue of $L$ on any compact domain of $M$
is positive (see \cite{FS, p.201}), we will use this fact from time to time. Then the
following comparison theorem holds.

\proclaim{Lemma 2.2} Assume that the first Dirichlet eigenvalue of $L$
on $M$ is nonnegative. Suppose $u$ and $v$ are two solutions of
(2.2) in a bounded domain   $\Omega$ such that $u\ge v$ on
$\p\Omega$. Then $u> v$ in $\Omega$ unless $u\equiv v$.

\endproclaim

\demo{Proof} Let $w=u-v$, then $w$ is also a solution of (2.2) and $w\ge 0$ on $\p\Omega$.
 Suppose $\inf_\Omega w<0$. Let $D=\{w<0\}$. Then $D$ is an open   subset of $\Omega$ and $w=0$ on $\p D$.
  We have

$$
\int_{D}\lf|\nabla w\ri|^2+q w^2=0.
$$
Since the first Dirichlet eigenvalue of $L$ in $M$ is nonnegative, the first eigenvalue of $L$ in $D$ must be
 positive. This is impossible. Hence $\min_\Omega w\ge0$.
 By the strong maximum principle \cite{GT, p. 35} we conclude that $w>0$ in $\Omega$ unless $w\equiv0$.
\enddemo

It is known that the first Dirichlet eigenvalue of $L$ is nonnegative if and only if (2.2) has a positive solution,
 see \cite{FS, p.201}. We want to discuss the asymptotically behavior of the positive solution when it exists.

\proclaim{Lemma 2.3} Let $q$ be a smooth function in $M$ such that $q=O(r^{-n})$ and let

$$
v(x)=-\int_M G(x,y)q(y)dy.
$$
Then $v$ is the unique solution of $\Delta v=q$ with $\lim_{x\to\infty} v(x)=0$. Moreover,
  $v=O(r^{2-n})$ and $\int_M|\nabla v|^2<\infty$.

\endproclaim

\demo{Proof} The fact that $v$ is the unique solution of $\Delta v=q$ with $\lim_{x\to\infty} v(x)=0$
      follows from the fact that $M$ is AF and Lemma 2.1. The fact that $v=O(r^{2-n})$ also follows from the fact
      that $M$ is AF. To prove that $|\nabla v|\in L^2(M)$, multiply $\Delta v=q$ by $\varphi^2 v$ and integrating
      by parts, where $\varphi$ is a cutoff function, we have

$$
\int_M\varphi^2|\nabla v|^2\le 4\int_M|\nabla
\varphi|^2v^2+\int_M|qv|\varphi^2.
$$
Using the fact that $v=O(r^{2-n})$, $q=0(r^{-n})$ it is easy to see that $\int_M|\nabla v|^2<\infty$.

\enddemo

\proclaim{Lemma 2.4} Suppose the first Dirichlet eigenvalue of $L$
on $M$ is nonnegative. Then the positive solution of (2.2)
 is unique in the sense that if $u$ and $v$ are two positive solutions of (2.2), then $u=\beta v$ for some $\beta>0$.
 Moreover, any positive solution of (2.2) is bounded and $\lim_{x\to\infty}u(x)$ exists. In fact, $u(x)=b+O(r^{2-n})$
 for some constant $b\ge0$.

\endproclaim

\demo{Proof} Let $u$ and $v$ be two positive solutions of (2.2). We may assume that there exist  $x_i\to\infty$
such that $u(x_i)\ge v(x_i)$. Suppose $r(x_i)=d(p,x_i)=R_i$ where $p$ is a fixed point. By Lemma 2.1, there is a
constant $C_1>0$ depending only on $M$ and $q$  such that

$$
\inf_{\p B_p(R_i)}u\ge C_1u(x_i)\ge C_1v(x_i)\ge C_1^2\sup_{\p
B_p(R_i)}v.
$$
By Lemma 2.2, we conclude that $u\ge C_1^2v$ in $B_p(R_i)$ for all $i$. Hence $u\ge C_1^2 v$ in $M$.
Let

$$
\beta^*=\sup\{\beta>0|\ u\ge \beta v \text{\ in $M$}\}.
$$
Then $u-\beta^* v$ is a solution of (2.2). Suppose $u\not\equiv \beta^*v$, then $w=u-\beta^* v>0$ in $M$ by Lemma 2.2.
Choose any $R'_i\to\infty$ and choose any $x'_i\in \p B_p(R'_i)$. Let
$a_i=\frac{w(x'_i)}{v(x'_i)}.$ We claim that $a_i$ is bounded from below away from zero. Otherwise,
passing to a subsequence, we may assume that $a_i\to0$. Since $w$ is also a positive solution of (2.2), we can
argue as before to conclude that

$$
\inf_{\p B_p(R'_i)}a_iv\ge C_1a_iv(x'_i)= C_1 w(x'_i)\ge
C_1^2\sup_{\p B_p(R'_i)}w.
$$
By Lemma 2.2, we have that $a_iv\ge C^2 w$ in $B_p(R'_i)$. From this it is easy to see that $w\equiv0$,
which is a contradiction. Hence $a_i\ge a>0$ for all $i$. Then as before

$$
\inf_{\p B_p(R'_i)}w\ge C_1w(x'_i)\ge C_1av(x'_i)\ge C_1^2\sup_{\p
B_p(R'_i)}v. 
$$
We conclude that $w\ge C_1^2av$ and $u\ge \lf(\beta^*+C_1^2a\ri)v$. This contradicts the definition of $\beta^*$.
So $u=\beta^*v$. This completes the proof of uniqueness.

It remains to prove the other assertions in the lemma.
Suppose $u$ is unbounded.   Then by Lemma 2.1, there exist
$r_i\to\infty$ such that $\inf_{\p B_p(r_i)}u\ge i.$ By Lemma 2.3,
let $\varphi$ be the solution of $\Delta \varphi=q$ such that $\varphi=O(r^{2-n})$.
Then

$$
\split \Delta \log u&=\frac{\Delta u}{u}-|\nabla \log u|^2\\&
=q-|\nabla \log u|^2\\&\le \Delta \varphi.
\endsplit
$$
Since $\inf_{\p B_p(r_i)}\log u\ge \log i$, $\varphi\le \log u-\log i+1$ in $B_p(r_i)$ for all $i$ provided $i$ is large. This is impossible. Hence $u$ is bounded.

Finally, $\Delta u=qu$ and $qu=O(r^{-n})$ because $u$ is bounded. Let $v'$ be the solution of $\Delta v'=qu$ given by
Lemma 2.3. Then $\Delta(u-v')=0$ and hence $u-v'$ must be constant by Lemma 2.1. Since $v'=O(r^{2-n})$, the lemma
follows.
\enddemo

Next, we will give a necessary and sufficient condition that the solution in Lemma 2.4 has a positive limit at infinity.

\proclaim{Theorem 2.1} Let $q$ be a smooth function on $M$ such that $q=O(r^{-n})$. Then

$$
\Delta u-qu=0
$$
has a positive solution $u$ satisfying $u=1+O(r^{2-n})$ if and only if there is a smooth function $f\ge0$,
$f\not\equiv0$ and $f=O(r^{-n})$ such that the operator $\Delta-q+f$ has nonnegative Dirichlet eigenvalue.

\endproclaim

\demo{Proof} Suppose there is a smooth function $f\ge0$, $f\not\equiv0$ and $f=O(r^{-n})$ such that the operator
$\Delta-q+f$ has nonnegative Dirichlet eigenvalue. Then the first Dirichlet eigenvalue of $\Delta-q$ is also nonnegative.
 By Lemma 2.4, we can find a positive solution of $Lu=0$ such that $ u=a+O(r^{2-n})$ for some $a\ge0$. Suppose $a=0$,
 then for any $0<\e<\max_M u$,   the set $\{u>\e\}$ is bounded in $M$ and the family of sets $\{u>\e\}$ with
 $\max_M u>\e>0$ exhausts $M$. Let $0<\e_0<\max_M u$ and let $R>0$ be fixed. Then for $0<\e\le \e_0$,
 $B_p(R)\subset \{u>\e\}$ if $\e$ is small enough, where $B_p(R)$ is the geodesic ball with center at a fixed point
 $p$   with radius $R$. For such $\e$, we have

$$
\split 0&\le
\int_{\{u>\e\}}\lf[|\nabla (u-\e)|^2+q(u-\e)^2 -f(u-\e)^2\ri]\\
&=-\e\int_{ \{u>\e\}}q(u-\e)-\int_{\{u>\e\}}f(u-\e)^2\\
&=- \e\int_{ \{u>\e\}}q u+\e^2\int_{ \{u>\e\}}q-\int_{\{u>\e\}}f(u-\e)^2\\
&=-\e\int_{ \{u>\e\}}q u+\e^2\int_{ \{u>\e\}}\frac{|\nabla u|^2}{u^2}+\e \int_{ \{u>\e\}}\frac{\p u}{\p\nu}-\int_{\{u>\e\}}f(u-\e)^2\\
&\le -\e\int_{ \{u>\e\}}q u+ \int_{ \{u>\e\}\setminus B_p(R)}
|\nabla u|^2 +\e^2\int_{ B_p(R)}\frac{|\nabla
u|^2}{u^2}-\int_{\{u>\e_0\}}f(u-\e)^2 \endsplit
$$
where $\nu$ is the unit outward normal of $\{u=\e\}$ so that $\frac{\p u}{\p \nu}\le 0$ and we have used the
fact that $f\ge0$. Let $\e\to0$ using the fact that $qu=O(r^{2-2n})$ so that $qu\in L^1(M)$ we have

$$
0\le \int_{ M\setminus B_o(R)} |\nabla u|^2
-\int_{\{u>\e_0\}}fu^2.
$$
By Lemma 2.3, $\int_M|\nabla u|^2<\infty$,  hence if we   let $R\to\infty$ and then let $\e_0\to0$, we have

$$
0\le -\int_M fu^2.
$$
Since $u>0$,  $f\ge0$ and $f\not\equiv0$, this is impossible.

Conversely, suppose there is a positive solution of $\Delta u-qu=
0$ with $u=1+O(r^{2-n})$. Let $q_+$ be the positive part of $q$,
then $q_+=O(r^{-n})$. Let $k$ be a smooth positive function such
that
 $k=0(r^{-n})$ and $k> q_+$. It is easy to see that the first Dirichlet eigenvalue of $\Delta -k$ is nonnegative.
 By Lemma 2.4, there is bounded and positive solution $w$ of $\Delta w-kw=0$. By multiplying $w$ by a positive number,
  we may assume that $ u-w\ge a>0$ for some $a$, where we have used the fact that $\inf_M u>0$ and $w$ is bounded.
  Let $f= ( \Delta w-qw)/(u-w)$. Since $k=O(r^{-n})$, $q=O(r^{-n})$, $w$ is bounded and $ u-w\ge a$,
  it is easy to see that $f=O(r^{-n})$. Also

$$
\Delta w-qw=kw-q_+w+q_-w>0
$$
because $k> q_+$, $w> 0$, where $q_-$ is the negative part of $q$. We conclude that   $f>0$. On the other hand,

$$
\split
\Delta(u-w)-(q-f)(u-w)&=\lf(\Delta u-qu\ri)-\Delta w+qw+f(u-w)\\
&=0.
\endsplit
$$
Hence the first Dirichlet eigenvalue of $\Delta-q+f$ is nonnegative by \cite{FS, p.201}.

\enddemo

\noindent{\bf Example}: It is easy to construct examples of $q$ such that $q=O(r^{-n})$ and the first
eigenvalue of $L=\Delta -q$ is nonnegative, but the positive solution $u$ of $Lu=0$ has the property that
 $\lim_{x\to\infty}u=0$. In fact, one can construct example with $q$ satisfying the additional property that $q$ has compact support and in particular $q\in L^1(M)$.
First, let $q'\le 0$ be a smooth function with compact support such that $q'<0$ somewhere. Let $u$ be a positive solution of
 $\Delta u=q'$ such that $u\to 0$ near infinity. Then $u$ is a positive solution of $\Delta u-q u=0$, where $q =q'/u$
 which is smooth with compact support. Note that by the uniqueness result in Lemma 2.4, for this $q$ every other positive
  solution of $\Delta v-qv=0$ must be asymptotically zero.\vskip .1cm

 Let $M$ be an AF manifold and let $\Lambda>0$ be the Sobolev constant on $M$. Namely,

$$
\Lambda=\inf\lf\{\frac{\int_M|\nabla
f|^2}{\lf(\int_M|f|^{\frac{2n}{n-2}}\ri)^{\frac {n}{n-2}}}\bigg|\
f\in C_0^\infty(M), \ f\not\equiv0\ri\}.\tag2.3
$$
Since $M$ is AF, $\Lambda>0$, Theorem 2.1 implies  the following result in \cite{SY1}.

\proclaim{Corollary 2.1} Let $(M^n,g)$ be an AF manifold and let $\Lambda$ be the Sobolev constant for $M$.
Let $q$ be a smooth function on $M$ such that $q=O(r^{-n})$. Suppose

$$
a= \lf(\int_Mq_-^{\frac{n}{2}}\ri)^{\frac{2}{n}} <\Lambda.$$
where $q_-$ is the negative part of $q$. Then $\Delta u-qu=0$ has
a positive solution $u$ such that $u=1+O(r^{2-n})$.

\endproclaim

\demo{Proof} For any smooth function $f$ with compact support,

$$
\split \int_M |\nabla f|^2+\int_M qf^2&=
\int_M |\nabla f|^2-\int_M q_- f^2+\int_Mq_+ f^2\\
&\ge \lf(\Lambda-a\ri)\lf(\int_M |f|^\frac {2n}{n-2}
\ri)^{\frac{n-2}n}.
\endsplit\tag2.4
$$
Since $\Lambda-a>0$, we can find a smooth function $h\ge0$ with compact support such that $h\not\equiv0$ and

$$
\lf(\int_Mh^{\frac{n}{2}}\ri)^{\frac{2}{n}}\le \Lambda-a.
$$

By (2.4), we have that

$$
\int_M |\nabla f|^2+\int_M qf^2\ge \lf(\int_Mh^{\frac{n}{2}}\ri)^{\frac{2}{n}}\lf(\int_M |f|^\frac {2n}{n-2}
\ri)^{\frac{n-2}n}\ge \int_Mhf^2,
$$
and hence the first Dirichlet eigenvalue of $\Delta-(q-h)$ is nonnegative. The corollary then follows from Theorem 2.1

\enddemo

\proclaim{Corollary 2.2} Let $(M^n,g)$ be an AF manifold with
scalar curvature $\cR$. Suppose the first Dirichlet eigenvalue of
$\Delta - \lf[(n-2)/4(n-1)\cR-f\ri]$ is nonnegative for some smooth
function $f\ge0$, $f\not\equiv0$ and $f=O(r^{-n})$. In addition,
suppose there is $\alpha>0$ such that the H\"older norm of $\cR$
with exponent $\alpha$ in $B_x(\frac 12r_x)$ decays like
$r_x^{-n-\alpha}$. Then $M$ is conformally scalar flat in the
sense that there is a smooth positive function $u$ such that
$(M,u^{\frac 4{n-2}}g)$ is an AF manifold with zero scalar
curvature.
\endproclaim
\demo{Proof} By Theorem 2.1 and the fact that $\cR=O(r^{-n})$,
there is a positive solution $u$ of $\Delta u-\cR u=0$ such that
$u(x)=1+O(r^{2-n})$  as $x\to\infty$. By the assumption on the
H\"older norm of $\cR$ and the fact that $g$ is AF, we have
$|\nabla u|=O(r^{1-n})$ and $|\nabla^2u|=O(r^{-n})$ by the
Schauder estimate \cite{GT, Theorem 6.2}. Hence  $u^{\frac
4{n-2}}g$ is also AF and has zero scalar curvature.
\enddemo

\proclaim{Remark 2.1} Suppose the first Dirichlet eigenvalue of
$\Delta -  (n-2)/4(n-1)\cR $ is nonnegative on $M$, then the first
Dirichlet eigenvalue of this operator will be positive on every
compact domain. Hence

$$
\int_M|\nabla f|^2+\frac{n-2}{4(n-1)}\cR f^2>0
$$
for all $f\in C_0^\infty(M)$. From the example before Corollary 2.1, one can see that this may not imply that
$M$ is conformally scalar flat. Hence the result in \cite{CB, Theorem 2.1, (I) implies (II)} seems to be incorrect. This   has also been noticed by Maxwell \cite{Ma}. In fact, Maxwell has obtained results similar to Corollary 2.2 with conditions   in terms of positivity of certain Sobolev quotient, see \cite{Proposition 4.1, Ma}. Moreover, AF manifolds with boundary are also discussed there and hence the results  are more general than Corollary 2.2.\endproclaim

Using  Theorem 2.1, one can find a harmonic spinor on $M$ which is asymptotically parallel near
infinity under certain condition on $\cR$. Denote the Dirac operator on $M$ by $\D$.

\proclaim{Theorem 2.2} Let $(M^n,g)$ be a spin AF manifold with scalar curvature $\cR$.
 Suppose there is a  smooth function $f\ge0$, $f\not\equiv0$ and $f=O(r^{-n})$ such that the first
 Dirichlet eigenvalue of $\Delta-(\frac14\cR-f)$ is nonnegative. Then for any spinor  $\Psi_0$ in $\R^n$
 which is parallel     with respect to the Euclidean metric, there exists a unique harmonic spinor $\Psi$ on
 $M$ such that $\D\Psi=0$ and there is a constant $C$ such that $|\Psi-\Psi_0|\le  Cr^{2-n}$ at infinity, where
 $|\cdot|$ is the norm with respect to $g$.

\endproclaim

\demo{Proof} Let $\Psi_0$ be a parallel spinor over $\R^n$ near
infinity.  $\Psi_0$ can be considered as a spinor over $M$ near
infinity. Extend $\Psi_0$ to be smooth on $M$. As in \cite{ST1,
\S3}, let $p$ be a fixed point, for any $R>0$ we first solve:

$$
\cases
 \D^2 \sigma_R &=-\D^2 \Psi_0 \text{\ \ in $B_p(R)$,}\\
\sigma_R|_{\p B_p(R)} &=0
\endcases\tag2.5
$$

To prove the solution exists, it is sufficient to prove that

$$a(\Psi,\Psi)=\int_{B(R)}\lgle \D\Psi,\D\Psi\rgle \geq \delta  \int_{B(R)}\lf(|\nabla \Psi|^2+|\Psi|^2\ri) \tag2.6
$$
for some $\delta>0$ for all $\Psi\in W^{1,2}_0(B_p(R))$. By Lichnerowicz formula:

$$
\split
a(\Psi,\Psi)&=\int_{B_p(R)}\lgle \nabla\Psi,\nabla\Psi\rgle +\frac14\cR\lgle \Psi,\Psi\rgle \\
& = \int_{B_p(R)}|\nabla \Psi|^2+\frac14\cR|\Psi|^2\\
&=(1-\tau)\lf[\int_{B_p(R)}|\nabla \Psi|^2+\frac1{4}\int_{B_p(R)}\cR|\Psi|^2\ri]+ \tau\lf[\int_{B_p(R)}|\nabla \Psi|^2+\frac1{4}\int_{B_p(R)}\cR|\Psi|^2\ri]\\
&\ge (1-\tau)\lf[\int_{B_p(R)}|\nabla \Psi|^2+\frac1{4}\int_{B_p(R)}\cR|\Psi|^2\ri]+\lambda\tau\int_{B(R)}|\Psi|^2\\
&\ge
(1-\tau)\int_{B(R)}|\nabla\Psi|^2+(\frac{1-\tau}{4}\inf_{B(R)}\cR+\lambda\tau)\int_{B(R)}
|\Psi|^2\endsplit
$$
where $0<\tau<1$ is a constant and $\lambda>0$ is the first eigenvalue of $\Delta -\frac14 \cR$ in $B_p(R)$ which is
 positive by assumption. Choose $\tau$ close enough to 1, (2.6) follows.

By Theorem 2.1, we can find a positive solution $u$ of

$$
\Delta u-\frac14\cR u=0\tag2.7
$$
such that $u\to1$ at infinity. Without loss of generality, we assume that the norm of $\Psi_0$ is 1 with respect to the Euclidean metric near infinity. Let $\Psi_R=\sigma_R+\Psi_0$, then $\D^2\Psi=0$  on $B_p(R)$. We want to prove that $\Psi_R$ is uniformly bounded. By Lichnerowicz formula, we have

$$
\frac12\Delta |\Psi_R|^2\ge\frac14\cR|\Psi_R|^2+|\nabla \Psi_R|^2
$$
and so

$$
 \Delta |\Psi_R|\ge\frac14\cR|\Psi_R|.
$$
Since $|\Psi_R|^2=|\Psi_0|^2$ on $\p B_p(R)$ and $|\Psi_0|$ is asymptotically  1 with respect to $g$ as $g$ is an AF metric, by Lemma 2.1,  we conclude that given any $\e>0$,

$$
|\Psi_R|\le u+\e
$$
in $B_p(R)$ if $R$ is large enough. Hence passing to a subsequence, $\Psi_R$ converges to a   spinor
$\Psi$ such that $\D^2\Psi=0$ and

$$
|\Psi|\le u\tag2.8
$$
on $M$.

To prove that $\Psi$ is harmonic, by Lemma 3.4 in \cite{ST1} it is sufficient to prove that

$$
\int_M|\D\Psi|^2<\infty.\tag2.9
$$
For each $R$, $\Psi_R-\Psi_0=\sigma_R=0$ on $\p B(R)$, so

$$
\int_{B(R)}\lgle \D\Psi_R,\D(\Psi_R-\Psi_0)\rgle=0.
$$
So

$$
\int_{B(R)}|\D\Psi_R|^2\le \int_{M}|\D\Psi_0|^2.
$$
Let $R\to\infty$, we see that (2.9) is true.

Next, we want to estimate $|\Psi-\Psi_0|$. By Lichnerowicz formula, for each $R$, we have

$$
\Delta|\Psi_R-\Psi_0|\ge -|\D^2\Psi_0|+\frac 14\cR|\Psi-\Psi_0|.
$$
Since $|\cR|=O(r^{-n})$ and $|\D^2\Psi_0|=O(r^{-n})$ and the fact that $|\Psi-\Psi_0|$ is bounded by (2.7),  we have that

$$
\Delta|\Psi_R-\Psi_0|\ge -C_1(1+r)^{-n}
$$
for some positive constant $C_1$. By Lemma 2.2, one can prove that   $|\Psi-\Psi_0|\le C_2r^{2-n}$ near infinity for some constant $C_2$.

Suppose $\Psi_1$ and $\Psi_2$ are two harmonic spinors such that $|\Psi_1-\Psi_0|\to0$ and $|\Psi_2-\Psi_0|\to0$ at infinity. By lichnerowicz formula, we have

$$
\Delta|\Psi_1-\Psi_2|\ge \frac14\cR |\Psi_1-\Psi_2|.
$$
By Lemma 2.1, we can conclude that $|\Psi_1-\Psi_2|\le \e u$ for any $\e>0$. Hence  $\Psi_1=\Psi_2$.
This completes the proof of the theorem.

\enddemo

For simplicity, the harmonic spinor obtained in the theorem is said to be the harmonic spinor with
 boundary value $\Psi_0$.

We can express the mass of $M$ in terms of the harmonic spinor as in  \cite{W, PT}:

\proclaim{Corollary 2.3} Let $(M^n,g)$ be an AF manifold satisfying the conditions of Theorem 2.2.
Let $\Psi$ be the harmonic spinor on $M$ with boundary value $\Psi_0$ where $\Psi_0$ is a parallel spinor
with respect to the Euclidean metric such that $|\Psi_0|=1$ near infinity. Then the mass of $M$ is given by

$$
c(n)\frak m=\int_M\lf(|\nabla \Psi|^2+\frac14\cR|\Psi|^2\ri).
$$
where $c(n)$ is a positive constant depending only on $n$.
\endproclaim
\subheading{\S3 Some lower bounds of ADM mass}\vskip .2cm

In this section, we shall give some lower bounds of the ADM mass of an AF manifold without assuming that
 the scalar curvature is nonnegative. Let $(M^n,g)$ be an AF manifold with one end as in \S2. In case
$n\ge 4$, we assume that $M$ is spin. Let $\Lambda>0$ be the
Sobolev constant on $M$ defined in (2.3).   Hence   if $f\in
C_0^\infty(M)$, then

$$
\int_M|\nabla f|^2\ge \Lambda
\lf(\int_M|f|^{\frac{2n}{n-2}}\ri)^{\frac {n}{n-2}}.\tag3.1
$$
It is easy to see that (3.2) is still true if $f$ is smooth such that $|f|=O(r^{-\tau})$ with $\tau>(n-2)/2$.

\proclaim{Theorem 3.1} Let $(M^n,g)$ be a spin AF manifold with scalar curvature $\cR$. Let $\cR_+$ and $\cR_-$ be the positive and negative part of $\cR$ respectively. Suppose

$$
a=\lf(\int_M\lf(\frac{\cR_-}4\ri)^{\frac{n}{2}}\ri)^{\frac{2}{n}}
<\frac\Lambda2.\tag3.2
$$
Then the mass $\frak m$ of $M$ has a lower bound given by

$$
\frak m\ge C(n)\frac{\Lambda -2a}{\Lambda
+b-a}\lf(\int_M\cR_+-\frac{\Lambda+2b}{\Lambda-2a}\int_M\cR_-\ri),
$$
where $C(n)$ is a positive constant depending only on $n$ and

$$
b=\lf(\int_M\lf(\frac{\cR_+}4\ri)^{\frac{n}{2}}\ri)^{\frac{2}{n}}.
$$
In particular, if   $\cR\ge 0$ then

$$
\frak m\ge C(n)\frac{\Lambda}{\Lambda +b}\int_M\cR.
$$
\endproclaim

\demo{Proof} By Corollary 2.1,  there is a positive solution $u$ of

$$
\Delta u-\frac{\cR}4u=0\tag3.3
$$
with $ u =1+(r^{2-n})$. By Theorem 2.2, there is a harmonic spinor $\Psi$ such that its norm is asymptotically equal to 1 near infinity. Namely, $|\Psi|=1+O(r^{2-n})$. Let $v=|\Psi|$, then the ADM mass of $M$ satisfies

$$
\split
C(n)\frak m&=\int_M\lf(|\nabla \Psi|^2+\frac{\cR}4|\Psi|^2\ri)\\
&\ge \int_M\lf(|\nabla v|^2+\frac{\cR}4v^2\ri)\\
&\ge  \int_M |\nabla (v-1) |^2+
 \lf(1-\frac1\e\ri)\int_M \frac{\cR_+}4  (v-1)^2-(1+\frac1\e)\int_M \frac{\cR_-}4(v-1)^2\\
&\quad+(1-\e)\int_M\frac{\cR_+}4 -(1+\e)\int_M\frac{\cR_-}4\\
&\ge \lf(\Lambda-\lf( \frac1\e-1\ri)b-(1+\frac1\e)a\ri)\lf(\int_M (v-1)^{\frac {2n}{n-2}}\ri)^{\frac {n-2}n}+(1-\e)\int_M\frac{\cR_+}4\\
&\qquad -(1+\e)\int_M\frac{\cR_-}4.
\endsplit\tag3.4
$$
for any  $0<\e\le1$, where we have used the fact that $|v-1|=O(r^{2-n})$ and the Sobolev inequality (3.1).
Choose $\e$ such that

$$
\Lambda-\lf( \frac1\e-1\ri)b-(1+\frac1\e)a=0.
$$
That is to say
$$
\e=\frac{b+a}{\Lambda+b-a}.
$$
Note that $0<\e\le 1$ unless $\cR\equiv0$ because $\Lambda>2a$. If $\cR\equiv0$, then the theorem is true by the positive mass theorem \cite{SY1-2,W,PT}. If $\cR\not\equiv0$, then the theorem follows from (3.4) and the definition of $\e$.
\enddemo

\proclaim{Corollary 3.1} Same assumptions and notations as in Theorem 3.1. Suppose
$$
\int_M\cR_+\ge\frac{\Lambda+2b}{\Lambda-2a}\int_M\cR_-
$$
then the mass $\frak m$ of $M$ is nonnegative. Moreover, $\frak m=0$ if and only if $M$ is the Euclidean space.
\endproclaim

\demo{Proof} The fact that $\frak m\ge0$ follows immediately from the theorem. Suppose $\frak m=0$,   then every inequality in (3.4) becomes an inequality. Hence  $\cR_+$, $(v-1)^2$ and $\cR_-$ are proportional to each other. So $\cR\equiv0$ and $M$ is the Euclidean space by the standard positive mass theorem.

\enddemo

\proclaim{Remark 3.1} The result of Corollary 3.1 can also be derived from the result in \cite{ZZ} under similar conditions. More precisely, if we replace $\cR_+/4$ and $\cR_-/4$ by $\cR_+/8$ and $\cR_-/8$ in the assumptions of the corollary, then the proof of the result  in \cite{ZZ, Theorem 4.1} together with similar derivation as in (3.4), we may also conclude that $\frak m\ge0$.
\endproclaim

Under the assumptions of Theorem 3.1 or the assumptions in \cite{ZZ, Theorem 4.1}, we must have $\int_M\cR\ge0$.
One might ask what might happen if $\int_M\cR<0$. In this situation, we  want to give a lower bound for the ADM mass
using the methods in \cite{BF, FK} and we shall give another condition so that the ADM mass is nonnegative. In the
following, we assume that $(M^n,g)$ is an AF manifold with scalar curvature $\cR$ which may be negative somewhere.
 We always assume that the operator $\Delta -1/4\cR+f$ has nonnegative Dirichlet eigenvalue  for some smooth function $
 f\ge0$, $f\not\equiv0$ and $f=O(r^{-n})$. Hence

$$
\Delta u-\frac14\cR u=0\tag3.5
$$
has a positive solution which tends to 1 near infinity by Theorem 2.1.  Denote

$$
A=\sup_M u.\tag3.6
$$

We will follow the arguments   in \cite{FK}.  Let $\Psi_0$ be a spinor which is parallel near infinity with respect to the Euclidean metric such that $|\Psi_0|_e=1$ near infinity, where $|\cdot|_e$ is the norm of $\Psi_0$ with respect to the Euclidean metric. By Theorem 2.2, because of (3.5), we can find a harmonic spinor $\Psi$ on $M$ such that $|\Psi-\Psi_0|\to 0$ and the ADM mass $\frak m$ is given by

$$
C(n)\frak m=\int_M\lf(|\nabla \Psi|^2+\frac14\cR|\psi|^2\ri).
$$
where $C(n)>0$ is a constant depending only on $n$. For such a $\Psi$, define

$$
\frak m_\Psi=C(n)\frak m-\frac14\int_M\cR|\Psi|^2=\int_M|\nabla
\Psi|^2.\tag3.7
$$
Hence $\frak m_\Psi\ge0$ for all such $\Psi$.
Let $\Psi_0$ and $\Psi$ as above. We have the following:

\proclaim{Lemma 3.1} For all $x\in M$,
$$
|\Psi_x|\le u(x)\le A \tag3.7
$$
where $u$ is the positive solution in (3.5) and $A$ is given by (3.6).
\endproclaim

\demo{Proof} By the Lichnerowicz formula,

$$
\frac 12\Delta|\Psi|^2=|\nabla \Psi|^2+\frac14\cR|\Psi|^2
$$
and so

$$
\Delta |\Psi|\ge \frac14\cR|\Psi|.
$$
 Since $|\Psi| \to 1$ near infinity, the lemma follows by maximum principle Lemma 2.2.
\enddemo

Let $\Psi$ be as in Lemma 3.1. We have:

\proclaim{Lemma 3.2} With the above notations, let $\eta\ge0$ be a smooth function on $M$ such that $\sup_M(|\eta|+|\Delta\eta|)<\infty,$ we have

$$
\int_M\eta||\nabla^2\Psi||^2\le C_1\frak m_\Psi \sup_M(|\eta
Rm|+|\Delta\eta|)+C_2A\sqrt{\frak m_\Psi}||\eta\nabla Rm||_2
\tag3.8
$$
for some constants $C_1$, $C_2$ depending only on $n$

\endproclaim

\demo{Proof} The proof is exactly as in Corollary 3.2 in \cite{FK}, except in the last part, we have to use Lemma 3.1 and the definition of $\frak m_\Psi$.
\enddemo

Let $N=2^{[\frac n2]}$ . Choose an orthonormal basis of constant spinors $\Psi_0^i$, $1\le i\le N$ with respect to the Euclidean metric and let $\Psi^i$ be the corresponding harmonic spinors. For $x\in M$, define $P_x$ as in \cite{FK}. Namely,

$$
P_x:S_x(M)\to S_x(M)
$$
where $S_x(M)$ is the fibre of the spinor bundle  associated with the spin structure through the spinor representation and

$$
P_x(\Psi)=\sum_{i=1}^n\langle
\Psi_x^i,\Psi_x\rangle\Psi^i_x.\tag3.9
$$
Note that $P_x$ will tend to the identity map near infinity.

\proclaim{Lemma 3.3} $|P_x|\le A $, where $|P_x|$ is the norm of the operator $P_x$ and $A$ is given by (3.6).
\endproclaim

\demo{Proof} The proof is same as Lemma 4.1 in \cite{FK}, except we have to use Lemma 3.1.
\enddemo

\proclaim{Lemma 3.4} There is a constant $c$ depending only on $n$ such that for any $\e>0$

$$
|| Id-P_x||^2<\e
$$
except on a set $D(\e)$ with

$$
\text{vol}(D(\e))\le\lf( \frac{c\sum_{i=1}^N\frak
m_{\Psi_i}}{\e^2\Lambda}\ri)^\frac{n}{n-2}
$$
where $Id$ is the identify map and $\Lambda$ is the Sobolev constant given in (3.1).

\endproclaim

\demo{Proof} The proof is exactly as in Lemma 4.2 in \cite{FK}.
\enddemo

Now choose $\e=N/32$, then outside $D=D(\e)$, $||Id-P||<\sqrt {N/32}$, and we have

$$
\frac N2|Rm|^2\le 32\sum_{i=1}^N||\nabla^2\Psi^i||^2
$$
by Lemma 5.1 in \cite{FK}. Combining this with (3.8), we have

$$
\split
\int_{M\setminus D)} \eta |Rm|^2&\le C(n)\int_{M\setminus D}\sum_{i=1}^N||\nabla^2\Psi^i||^2\\
&\le C_1\lf(\sum_{i=1}^N\frak m_{\Psi^i}\ri) \sup_M(|\eta
Rm|+|\Delta\eta|)+C_2A\sum_{i=1}^N\sqrt{\frak
m_{\Psi^i}}||\eta\nabla Rm||_2
\endsplit $$
where $C_1$ and $C_2$ are constants depending only on $n$.  Hence we have the following:

\proclaim{Theorem 3.2} Let $\Psi^i$, $1\le i\le N$ as above. Then there are constants $C_1(n)$, $C_2(n)$ depending only on $n$ such that for any smooth function $\eta$ on $M$ with
$\sup_M\lf(|\eta|+|\Delta \eta|\ri)<\infty$, we have that

$$
\int_{M\setminus D } \eta |Rm|^2 \le C_1\lf(\sum_{i=1}^N\frak
m_{\Psi^i}\ri) \sup_M(|\eta
Rm|+|\Delta\eta|)+C_2A\sum_{i=1}^N\sqrt{\frak
m_{\Psi^i}}||\eta\nabla Rm||_2.\tag3.10
$$
where $ D$ is a  set with

$$
\text{vol}(D)\le\lf( \frac{c\sum_{i=1}^N\frak
m_{\Psi_i}}{\lf(\frac{N}{32}\ri)^2 \Lambda}\ri)^\frac{n}{n-2}
$$
where $c$ is the constant in Lemma 3.4.
\endproclaim

Let  $\cR_-$ is the negative part of $\cR$ and let

$$
B=\inf_{D}\lf\{\int_{M\setminus D}|Rm|^2\bigg|\ \ \text{vol}(D)\le
\lf[\frac{cAN\int_M\cR_-}{\lf(\frac{N}{32}\ri)^2\Lambda}\ri]^{\frac{n}{n-2}}\ri\}.\tag3.11
$$

\proclaim{Corollary 3.2}  Let $(M^n,g)$ be a spin AF manifold with scalar curvature $\cR$ such that the operator $\Delta-\frac14\cR+f$ has nonnegative Dirichlet eigenvalue  for some smooth function $f\ge0$, $f\not\equiv0$, $f=O(r^{-n})$.
There exists $C(n)>0$ depending only on $n$ such that if

$$
B\ge C(n) A^2\lf[\sup_M|Rm|\int_MR_-+||\nabla Rm||_2\lf(
\int_MR_-\ri)^\frac12  \ri]
$$
then the mass of $M$ is nonnegative, where $A$ is given by (3.6) and  $B$ is given by (3.11).
\endproclaim

\demo{Proof} If $M$ is flat, then it is obvious that $\frak m=0$.
Suppose $M$ is non flat and suppose $\frak m<0$, then $\cR<0$
somewhere by the positive mass theorem. From the definitions of
$\frak m_{\Psi^i}$ and $A$ and by Lemma 3.1, we have that

$$
\frak m_{\Psi^i}<\frac {A^2}4\int_MR_- \tag3.12
$$
for all $1\le i\le N$.
Take $\eta\equiv1$ in (3.10), we have that

$$
\split
B&\le \int_{M\setminus D}|Rm|^2\\
&\le C_1\lf[\lf(\sum_{i=1}^N\frak m_{\Psi^i}\ri)\sup_M|Rm|+A\lf(\sum_{i=1}^N\sqrt{\frak m_{\psi^i}}\ri)  ||\nabla Rm||_2\ri]\\
&<  C_2 A^2\lf[\sup_M|Rm|\int_MR_-+||\nabla Rm||_2\lf(
\int_MR_-\ri)^\frac12  \ri].
\endsplit
$$
for some positive constants $C_1$, $C_2$ depending only on $n$ because $M$ is nonflat. 
From this, the result follows.

\enddemo

Under certain conditions, we can estimate $A$ from above. For example, we have the following:

\proclaim{Corollary 3.3} Let $(M^3,g)$ be an AF manifold with scalar curvature $\cR$ such that

$$\lf(\int_M\lf(\frac{\cR_-}{4}\ri)^\frac32\ri)^\frac23<\Lambda,$$
where $\Lambda$ is the Sobolev constant of $M$. Then there is a constant $C>0$ depending only on $n$ and positive constant $A$ depending only on $\cR_-$ and $\Lambda$ such that if
$$
B\ge CA^2 \lf[\sup_M|Rm|\int_MR_-+||\nabla Rm||_2\lf(
\int_MR_-\ri)^\frac12  \ri]
$$
where $B$ is as in Corollary 3.2.
Then $\frak m\ge0$.
\endproclaim
\demo{Proof} By Corollary 2.1, $M$ satisfies the conditions in Corollary 3.2.
It remains to prove that $A$ in Corollary 3.2 is less than some constant $C$ depending only on
$\cR_-$ and $\Lambda$. This will be proved in Corollary 4.1, next section.
\enddemo
\subheading{\S4 Nonnegativity and   lower estimates of Brown-York mass}\vskip .2cm

Suppose $(\Omega^3,g)$ is a $3$-dimensional compact manifold with
smooth boundary. For simplicity, in this section we always assume
that $\p\Omega=\Sigma$ is connected. Suppose the   Gauss curvature
of $\Sigma$ is positive, then by a classical result we know that
$\Sigma$ can be isometrically embedding into $\Bbb{R}^3$, see
\cite{N}.  Let $H_0$ be the mean curvature  of the embedding image
in $\Bbb{R}^3$ with respect to the  outward unit norm. Then the
Brown-York  mass of $(\Omega,g)$ is defined as follow (see
\cite{BY}):

$$\BY(\Omega)= \int_\Sigma (H_0 -H)d\sigma,\tag4.1$$
where $d\sigma$ is the volume element of $\Sigma$, $H$ is the mean
curvature of $\p \Omega$ with respect to original metric $g$ and
outward unit normal.

In \cite{ST1}, it was proved  that if (i) the scalar curvature of $(\Omega,g)$
is nonnegative and (ii) $H>0$,  then $\BY(\Omega) \geq 0$.
Moreover, equality holds if and only if $\Omega$ is a domain in
$\R^3$. In this section, we will discuss cases that either (i) or (ii) does not hold. As for ADM mass of an AF manifold,
we will also give some lower estimates of the Brown-York mass of a bounded domain.

In \cite{ST2}, some results in \cite{ST1} are generalized to the case that the Gauss curvature of $\Sigma$ is
only assumed to be nonnegative.  We believe that some of the results in this section are still true if $\Sigma$ is
only assumed to have nonnegative Gauss curvature. However, we always assume that the Gauss curvature of $\Sigma$ is
positive for simplicity.

Let us first consider the case that the scalar curvature $\cR$ may be negative somewhere. The idea is to solve

$$\cases
\Delta u-qu&=0\ \text{in $\Omega$}\\
u&=1\ \text{on $\p\Omega$}
\endcases\tag4.2
$$
If (4.2) has a positive solution with $q=\frac\cR8$, then the
metric $g_1=u^4g$ has zero scalar curvature. We can then apply the
result of \cite{ST1} provided that the mean curvature of $\Sigma$
with respect to $g_1$ is positive. This will be true if
$\frac{\partial u}{\partial\nu}$ is not too negative, where $\nu$
is the unit outward normal of $\Sigma$.
 To this end, we first get some upper estimates of the solution $u$ of (4.2). We have the following:

\proclaim{Lemma 4.1} Let $(\Omega^3,g)$ be a compact manifold with smooth boundary and let $q$ be a smooth
 function defined on $\ol \Omega$. Suppose
$$
\beta<\Lambda\tag4.3
$$
then (4.2) has a unique solution $u$ and such that
$$
0< u\le   1+ 27^\frac18\gamma\lf[\frac{(\alpha+1)(1+\Lambda-\beta)}{\Lambda(
\Lambda-\beta)}+1\ri]
 \tag4.4
$$
in $\Omega$ where $\alpha=\max_\Omega q_-$,
$\beta=\lf(\int_\Omega q_-^\frac32\ri)^\frac23$,
$
\gamma=\sup_{p\ge 1}\lf(\int_\Omega q_-^p\ri)^\frac1p
$, $q_-$ is the negative part of $q$ and $\Lambda$ is the Sobolev constant of $\Omega$.
 \endproclaim

\demo{Proof} Let $f\in C_0^\infty(\Omega)$ and $f\not\equiv 0$, then

$$
\split
\int_\Omega\lf(|\nabla f|^2+qf^2\ri)&\ge \int_\Omega\lf(|\nabla f|^2-q_-f^2\ri)\\
&\ge \lf(\Lambda-\beta\ri)\lf(\int_\Omega f^6\ri)^\frac13\\
&\ge \lf(\Lambda-\beta\ri)|\Omega|^{-\frac23}\int_\Omega f^2.
\endsplit
$$
Since $ \Lambda-\beta>0$,   (4.2) has a unique positive solution $u$ by \cite{GT, Theorem?}.
Let $v=u-1$, then
$$
\cases
\Delta v-qv&=q\ \text{in $\Omega$}\\
v&=0\ \text{on $\p\Omega$}.
\endcases\tag4.5
$$
Let $v_+=\max\{v,0\}$. For any integer $k\ge 1$, multiply the first equation in (4.5) by $v_+^{2k-1}$ and integrating by parts, we have that

$$
\split \int_\Omega q_-\lf(v_+^{2k-1}+v_+^{2k}\ri)
&\ge - \int_\Omega q \lf(v_+^{2k-1}+v_+^{2k}\ri)\\
&=(2k-1)\int_\Omega v_+^{2k-2}|\nabla v_+|^2\\
&=\frac{2k-1}{k^2}\int_\Omega|\nabla v_+^k|^2\\
&\ge \frac{\Lambda}{k }\lf(\int_\Omega v_+^{6k}\ri)^\frac13
\endsplit
$$
where we have used the fact that $v=0$ on $\p\Omega$. Hence

$$
\split
\lf(\int_\Omega v_+^{6k}\ri)^\frac13&\le \frac{k}{\Lambda}\lf[\alpha\int_\Omega v_+^{2k}+\lf(\int_\Omega v_+^{2k}\ri)^{\frac{2k-1}{2k}}\lf(\int_\Omega q_-^{2k}\ri)^{\frac1{2k}}\ri]\\
&\le \frac{k}{\Lambda}\lf[(\frac{2k-1}{2k}+\alpha)\int_\Omega v_+^{2k}+\frac1{2k}\int_\Omega q_-^{2k}\ri]\\
&\le \frac{(\alpha+1)k}{\Lambda}\max\lf\{\int_\Omega
v_+^{2k},\int_\Omega q_-^{2k}\ri\}
\endsplit\tag4.6
$$
where we have used H\"older and Young inequalities. For any $k\ge1$, let

$$
I_k=\lf( \int_\Omega v_+^{2k}\ri)^\frac1{2k}
$$
and let $a=(\alpha+1)/\Lambda$.
By (4.6) and the definition of $\gamma$, we have that

$$
\split
I_{3k}&\le \lf(\frac{(\alpha+1)k}{\Lambda}\ri)^\frac{1}{2k} \max\{I_k,\gamma\} \\
&=\lf(ak\ri)^{\frac1{2k}}\max\{I_k,\gamma\}.
\endsplit\tag4.7
$$
Let $\ell_0\ge0$ be such that $ak_0=a3^{\ell_0}\ge1$.   We claim that for $\ell\ge1$,

$$
\split
I_{3^\ell k_0}&\le a^{\frac{1}{2k_0}\sum_{i=0}^{\ell-1}3^{-i}}k_0^{\frac{1}{2k_0}\sum_{i=0}^{\ell-1}3^{-i}}3^{\frac{1}{2k_0}\sum_{i=0}^{\ell-1}i3^{-i}}\max\{I_{k_0},\gamma\}\\
&=a^{\frac{1}{2}\sum_{i=\ell_0}^{\ell_0+\ell-1}3^{-i}}3^{\frac12\sum_{i=\ell_0}^{\ell_0+\ell-1}i3^{-i}}\max\{I_{k_0},\gamma\}.
\endsplit\tag4.8
$$
 If $\ell=1$, then (4.8) is true by (4.7). Suppose (4.8) is true for $\ell$,   by (4.7), we have that

$$
\split
I&_{3^{\ell+1}k_0}\\
&\le \lf(a\cdot 3^\ell k_0 \ri)^\frac{1}{2\cdot 3^\ell k_0}\max\{I_{3^\ell k_0},\gamma\}\\
&\le \lf(a\cdot 3^\ell k_0 \ri)^\frac{1}{2\cdot 3^\ell k_0} \max\lf\{a^{\frac{1}{2k_0}\sum_{i=0}^{\ell-1}3^{-i}}k_0^{\frac{1}{2k_0}\sum_{i=0}^{\ell-1}3^{-i}}3^{\frac{1}{2k_0}\sum_{i=0}^{\ell-1}i3^{-i}}\max\{I_{k_0},\gamma\},\gamma\ri\}\\
&\le
a^{\frac{1}{2k_0}\sum_{i=0}^{\ell}3^{-i}}k_0^{\frac{1}{2k_0}\sum_{i=0}^{\ell}3^{-i}}3^{\frac{1}{2k_0}\sum_{i=0}^{\ell}i3^{-i}}\max\{I_{k_0},\gamma\}
\endsplit\tag4.9
$$
where we have used the fact that $ak_0\ge1$. Hence (4.8) is true. Suppose $\ell_0$ satisfies $a3^\ell<1$ for all $0\le \ell <\ell_0$, then for $k_0=3^{\ell_0}$,   (4.7) implies

$$
\split
I_{k_0}&=I_{3^{\ell_0}}\\
&\le \lf(a\cdot 3^{\ell_0-1} \ri)^\frac{1}{2\cdot 3^{\ell_0-1}}\max\{ I_{3^{\ell_0-1}},\gamma\}\\
&\le \max\lf\{\lf(a\cdot 3^{\ell_0-1} \ri)^\frac{1}{2\cdot 3^{\ell_0-1}}  I_{3^{\ell_0-1}},\gamma\ri\} \\
&\le \cdots\cdots\\
&\le \max\lf\{ a^{\frac{1}{2}\sum_{i=0}^{\ell_0-1}3^{-i}}
3^{\frac{1}{2}\sum_{i=1}^{\ell_0-1}i3^{-i}} I_3,\gamma \ri\}
\endsplit\tag4.10
$$
where we have used the fact that $a\cdot 3^\ell\le1$ if $\ell<\ell_0$. Note that

$$
\lim_{\ell\to\infty}I_{3^\ell}=\lim_{\ell\to\infty}\lf(
\int_\Omega
v_+^{3^\ell}\ri)^{\frac1{3^\ell}}=\lim_{\ell\to\infty}\lf(\frac{1}{|\Omega|}\int_\Omega
v_+^{3^\ell}\ri)^{\frac1{3^\ell}}=\sup_\Omega v^+.\tag4.11
$$
Suppose $a\ge1$,  we take $\ell_0=1$ in (4.8), we have

$$
I_{3^{\ell}}\le a^{\frac{1}{2}\sum_{i=1}^{ \ell
}3^{-i}}3^{\frac12\sum_{i=1}^{ \ell }i3^{-i}}\max\{ I_{3}, \gamma\}.
$$
Let $\ell\to\infty$, we have

$$
\sup_\Omega v_+\le 27^\frac18 a\max\{I_3,\gamma\}\le
27^{\frac18} (aI_3+(a+1)\gamma)\tag4.12
$$
because $3^{\sum_{i=1}^{ \infty }i3^{-i}}=27^\frac14$. Suppose $a<1$. Then choose $\ell_0\ge0$ such that $ak_0=a3^{\ell_0}\ge1$ for $\ell\ge\ell_0$ and such that $a3^\ell<1$ for $0\le \ell<\ell_0$. 
By (4.8) and (4.10), we see that (4.12) is still true.

To estimate $I_3$, as before we have
$$
\split
\Lambda\lf(\int_\Omega v_+^6\ri)^\frac13&\le \int_{\Omega}q_-\lf(v_+^2+v_+\ri)\\
&\le \lf(\int_\Omega q_-^\frac32\ri)^\frac23\lf(\int_\Omega
v_+^6\ri)^\frac13+\lf(\int_\Omega
q_-^\frac65\ri)^\frac56\lf(\int_\Omega v_+^6\ri)^\frac16
\endsplit
$$
By the defintions of $\beta$ and $\gamma$, we have
$$
(\Lambda-\beta)I_3\le \gamma.
$$
Combining this with (4.12), we have

$$
\sup_\Omega v_+\le
27^\frac18\gamma\lf[\frac{(\alpha+1)(1+\Lambda-\beta)}{\Lambda(
\Lambda-\beta)}+1\ri].$$
\enddemo

\proclaim{Corollary 4.1} Let $(M^3,g)$ be an AF manifold with Sobolev constant $\Lambda$. Let $q$ be a smooth   function such that $q=O(r^{-3})$,   and suppose

$$
\beta=\lf(\int_M q_-^\frac32\ri)^\frac23<\Lambda
$$
Then $\Delta -q$ has a positive solution $u$ which is asymptotically 1 near infinity such that

$$
0<u\le 1+ 27^\frac18\gamma\lf[\frac{(\alpha+1)(1+\Lambda-\beta)}{\Lambda(
\Lambda-\beta)}+1\ri]
$$
where $\alpha=\sup_M q_-$ and   $\gamma=\sup_{p\ge
1}\lf(\int_M q_-^p\ri)^\frac1p. $ In particular, $u$ is
bounded from above by a constant depending only on $q_-$ and $\Lambda$.

\endproclaim

\demo{Proof} The existence of $u$ which is asymptotically 1 near infinity is a consequence of Corollary 2.1. It is easy to see that $u$ is the limit of a sequence of solutions of

$$\cases
\Delta u_k-qu_k&=0\ \text{in $\Omega_k$}\\
u_k&=1\ \text{on $\p\Omega_k$}
\endcases
$$
where $\{\Omega_k\}_{k=1}^\infty$ is some family of bounded domains with smooth boundaries which exhausts $M$. The estimate of $u$ follows from the lemma.
\enddemo
\proclaim{Lemma 4.2} Let $(\Omega^3,g)$ be a compact manifold with smooth boundary $\Sigma$ and with scalar curvature $\cR$. Let
$$
I=\inf\lf\{\int_\Omega\lf(|\nabla w|^2+\frac \cR8w^2\ri)\bigg|\
\text{$w$ is smooth in\ $\ol\Omega$, $w\equiv1$\  on
$\p\Omega$}\ri\}.
$$

\roster
\item"{(i)}" Suppose  $\beta=\lf(\int_\Omega q_-^{\frac 32}\ri)^\frac23<  \frac{\Lambda}2 $,
 where $\Lambda$ is the Sobolev constant of $\Omega$, $q=\cR/8$.
Then
$$
I\ge\frac{\Lambda -2\beta}{8(\Lambda
+\delta-\beta)}\lf(\int_\Omega
\cR_+-\frac{\Lambda+2\delta}{\Lambda-2\beta}\int_\Omega\cR_-\ri)
$$
where $ \delta=\lf(\int_\Omega q_+^{\frac 32}\ri)^\frac23  $.
\item"{(ii)}" Let $\lambda$ be the first Dirichlet eigenvalue for the Laplacian of $\Omega$. Suppose $8\lambda+\cR>0$ in $\Omega$. Then
$$
I\ge \lambda\int_\Omega \frac{\cR}{8\lambda+\cR}.
$$
\endroster
\endproclaim
\demo{Proof} The proof of part (i) is similar to the proof of Theorem 3.1.

To prove (ii), let $w$ be smooth in $\ol\Omega$, $w\equiv1$ on $\p\Omega$, and let $v=w-1$, then

$$
\split
\int_\Omega\lf(|\nabla w|^2+\frac\cR8 w^2\ri)&=\int_\Omega\lf(|\nabla v|^2+\frac\cR8(1+v)^2\ri)\\
&\ge \frac18\int_\Omega\lf(8{\lambda}v^2+\cR(1+v)^2\ri)\\
&\ge \int_\Omega\frac{\lambda R}{8\lambda+\cR } .
\endsplit
$$
by minimizing $f(v)=8\lambda v^2+\cR (1+v)^2$, where we have used the fact that $8\lambda +\cR>0$ in $\Omega$. From this the   lemma follows.

\enddemo

Now consider a compact manifold $(\Omega,g)$ with smooth boundary
$\Sigma$  with positive mean curvature with respect to the outward
normal. Let $s_0>0$ be such that $d(x,\p\Omega)$ is smooth in
$\Omega_{s_0}=\{x|\ 0<d(x,\p\Omega)\le s_0\}$
 and such that the mean curvature  of $\{x|\ d(x,\p\Omega)=s\}$ with respect to the outward norm is positive for all $0\le s\le s_0$. We have the following theorem which implies that if $\cR$ is not too negative, then $\BY(\Omega)$ is still nonnegative.

\proclaim{Theorem 4.1} Let $(\Omega^3,g)$ be a compact manifold with boundary with scalar curvature $\cR$. With the above notations and assumptions, let $H_{\min}$ be the minimum of the mean curvature of $\Sigma$ with respect to the outward normal which is assumed to be positive and let

$$
\xi=\min\{\frac{\pi}{6s_0}, \frac12 H_{\min}\}.\tag4.13
$$
Suppose the Gauss curvature of $\Sigma$ is positive and suppose the following are true.

\roster \item"{(i)}" $\beta=\lf(\int_\Omega q_-^{\frac
32}\ri)^\frac23<  {\Lambda} $, where $\Lambda$ is the Sobolev
constant of $\Omega$ and $q=\cR/8$. 
\item"{(ii)}" $\alpha=\max_\Omega q_-\le \xi^2$. \item"{(iii)}"
$\gamma=\sup_{p\ge 1}\lf(\int_\Omega q_-^p\ri)^\frac1p
 \le  27^{-\frac18}\lf[\frac{(\alpha+1)(1+\Lambda-\beta)}{\Lambda-\beta}+1\ri]^{-1}\cdot\frac{\xi s_0}{10}.$
\endroster

Then the Brown-York mass of $\Omega$ satisfies:

$$
\BY(\Omega)\ge \frac14\inf\lf\{\int_\Omega\lf(|\nabla w|^2+\frac
\cR8w^2\ri)\bigg|\ w\equiv1\text{\ on $\Sigma$}\ri\}.
$$
\endproclaim

\proclaim{Remark 4.1} (a) The conditions of the theorem are
obviously satisfied if $\Omega$ has nonnegative scalar curvature.
(b) The conditions  (i)--(iii)  of the theorem will be  satisfied if $\alpha\le \xi^2$,   $\alpha|\Omega|^\frac23\le \Lambda/2$ and
 $\alpha\max\{|\Omega|,1\}\le 27^{-\frac18}\lf[\frac
{(\xi^2+1)(2+\Lambda)}{\Lambda^2}+1\ri]^{-1}\cdot\frac{\xi s_0}{10}.$  Hence,   there is a constant $C$ depending only on $\Lambda$,
$s_0$, $H_{\min}$ and $|\Omega|$ such that if $\alpha\le C$ then the
conditions of the theorem will be satisfied.

\endproclaim

\demo{Proof of Theorem 4.1} By (i) and Lemma 4.1, there is a
unique positive solution $u$ of (4.2) with $q=\cR/8$. Let
$\phi(s)=\cos\xi s+\sin \xi s$, $0\le s\le s_0$. On $0\le s\le
s_0$, $\phi>0$ by (4.13),
$$
\frac{d^2\phi}{ds^2}=-\xi^2\phi,\tag4.14
$$
and
$$
\split
\frac{d\phi}{ds}&=\xi\lf(-\sin\xi s+\cos\xi s\ri)\\
&\ge \xi \lf(-\frac{1}2+\frac{\sqrt 3}2\ri)\\
&\ge \frac{\xi}{10}.
\endsplit\tag4.15
$$
Moreover,
$$
\phi(0)=1,\ \phi(s_0)\ge 1+ \frac{\xi s_0}{10}, \
\phi'(0)=\xi\tag4.16
$$
Define $f(x)=\phi\lf(d(x,\p\Omega)\ri)$ for $x\in \Omega_{s_0}$. Then in $\Omega_{s_0}$

$$
\split
\Delta f-qf&=\phi''-H\phi'-q \phi\\
&\le \lf(-\xi^2+q_-\ri)\phi\\
&\le 0
\endsplit\tag4.17
$$
where $H$ is the mean curvature of $\{d(x,\p\Omega)=s\}$ and  we have used the facts that $H>0$, $\phi'\ge0$, $\phi\ge0$,   and $a=\max_\Omega q_-\le \xi^2$.
Moreover, $f=1$  on $\p\Omega$, and on $\{d(x,\p\Omega)=s_0\}$,

$$
f\ge 1+\frac{\xi s_0}{10}\ge u
$$
by (4.15), condition (iii) and Lemma 4.1. Since the first Dirichlet eigenvalue of $\Delta -q$ in $\Omega_{s_0}$ is positive, by the maximum principle, we have $f\ge u$ in $\Omega_{s_0}$. Hence on $\p\Omega$
$$
\frac{\p u}{\p\nu}\ge \frac{\p f}{\p\nu}=-\xi 
$$
where $\nu$ is  the unit outward normal of $\Sigma$.
Consider the metric $g_1=u^4g$, then the scalar curvature of $g_1$ is zero and the mean curvature $\ol H$ with respect to $g_1$ of $\p\Omega$ is
$$
\ol H=H+\frac{\p u}{\p\nu}\ge H-\xi>0
$$
where $H$ is the mean curvature with respect to  $g$ and we have used (4.13).
Since $u=1$ on the boundary, the induced metric on $\Sigma$ is the same as before and so the Gauss curvature of $\Sigma$ is positive. By \cite{ST1}, we have

$$
\int_{\Sigma}\lf(H_0-\ol H\ri)d\sigma\ge0.
$$
Since $\ol H=H+\frac14\frac{\p u}{\p\nu}$
$$
\BY(\Omega)=\int_{\Sigma}\lf(H_0-H\ri)d\sigma\ge
\frac14\int_{\Sigma}\frac{\p u}{\p\nu}d\sigma=\frac14\int_\Omega
|\nabla u|^2+qu^2. 
$$
From this the result follows.
\enddemo

We should remark that in Lemma 4.1 without assuming
$\beta<\Lambda$, we may obtained an upper bound for $u$ in terms
of $\Lambda$, $\alpha$, $|\Omega|$ and the first Dirichlet
eigenvalue of $\Delta$ provided   (4.2) has a positive solution
with $q=\cR/8$. Hence we may have a result similar to Theorem 4.1
without assuming $\beta<\Lambda$, provided (4.2) has a positive
solution with $q=\cR/8$.

By Theorem 4.1 and Lemma 4.2, we have:
\proclaim{Corollary 4.2} With the same assumptions and notations as in Theorem 4.1. 
\roster
\item"{(a)}" Suppose $\beta=\lf(\int_\Omega\lf(\frac{\cR_-}8\ri)^\frac32\ri)^\frac23<\frac12\Lambda$, then
$$
\BY(\Omega)\ge\frac{\Lambda -2\beta}{32(\Lambda
+\delta-\beta)}\lf(\int_\Omega
\cR_+-\frac{\Lambda+2\delta}{\Lambda-2\beta}\int_\Omega\cR_-\ri) \tag4.18
$$
where $\delta=\lf(\int_\Omega \lf(\frac{\cR_+}{8}\ri)^{\frac
32}\ri)^\frac23$.
In particular, $\BY(\Omega)\ge0$ if
$$\int_\Omega \cR_+\ge\frac{\Lambda+2\delta}{\Lambda-2\beta}\int_\Omega\cR_-.
$$

\item"{(b)}" Suppose    $8\lambda+\cR>0$ where $\lambda$ is the first Dirichlet eigenvalue for the Laplacian of $\Omega$, then 

$$
\BY(\Omega)\ge\frac{\lambda}4\int_\Omega
\frac{\cR}{8\lambda+\cR}.\tag4.19
$$
In particular, if $\int_\Omega \frac{\cR}{8\lambda+\cR}\ge0$, then $\BY(\Omega)\ge0$.
\endroster
\endproclaim

Hence if the mean curvature of $\Sigma$ is positive, then the Brown-York mass of $\Omega$ is still nonnegative if  the scalar curvature of $\Omega$ is not too negative.

Next we consider the case that the scalar curvature is positive but   the mean curvature $H$ of $\Sigma$ is negative  somewhere. We want to prove that if $H$ is not too negative then   $\BY(\Omega)$ is still nonnegative.  Let us fisrt consider a special case that $\cR\ge0$, but the mean curvature of the boundary is only assumed to be nonnegative, we have a simple proof of the nonnegativity of the Brown-York mass.

\proclaim{Proposition 4.1}
Let $(\Omega,g)$ be a $3$-dimensional compact Riemannian manifold
with smooth boundary $\Sigma$ such that $\Sigma$ has positive
Gauss curvature  and nonnegative   the mean curvature is only
assumed to satisfy $H\ge0$. Then
$$
\BY(\Omega)\ge0.
$$
In fact $\BY(\Omega)$ can be bounded from below by (4.18) and (4.19) with $\cR_-=0$.
\endproclaim

\demo{Proof}    Suppose $H>0$ at some point $p\in \Sigma$. Let $U$
be an open neighborhood of $p$ in $\Sigma$ such that $H\ge a>0$ in
$U$ for some positive constant $a>0$. Let $1\ge\varphi\ge0$ be
smooth function on $\Sigma$ with compact support in $U$ such that
$\varphi(p)>0$. For any $\epsilon>0$, Let $v_\epsilon$ be the
solution of
$$
\cases \Delta v_\epsilon-\frac R8v_\e&=0 \text{  in $\Omega$}\\
v_\epsilon &=1-\epsilon\varphi\text{ on $\Sigma$}.
\endcases
$$
Then $1>v_\epsilon>1-\epsilon$  by maximum principle. Let
$g_\epsilon=v_\epsilon^4g$. Then $(M,g_\epsilon)$ has zero scalar
curvature, such that the mean curvature of $\Sigma$ is positive
provided $\epsilon>0$ is small enough, because $v_\epsilon\to 1$
in $C^2$ and $v_\epsilon$ is not constant. The Gauss curvature of
$\Sigma$ is also positive. Hence

$$
\int_\Sigma \lf(H_0^\epsilon-H^\epsilon\ri)\,d\sigma_\e\ge0
$$
where $H_0^\e$, $H^\e$ are the mean curvatures of $\Sigma$ when embedded in $\R^3$ and in $\Omega$ respectively, $d\sigma_\e$ is the area element of with respect to $g_\e$. Let $\epsilon\to0$ and use Corollary 4.2, we conclude the proposition is true.

Suppose $H\equiv0$. If $\cR\equiv0$, then it is obvious the proposition is true. If $\cR>0$ somewhere, then the solution of

$$
\cases \Delta u-\frac R8u&=0 \text{  in $\Omega$}\\
u &=1 \text{ on $\Sigma$}.
\endcases
$$
satisfies $\frac{\p u}{\p \nu}>0$ by the strong maximum principle. In this case, the result follows as in Corollary 4.2.
\enddemo

For a more general case, as before let $s_0>0$ be such that $d(x,\p\Omega)$ is smooth in $\Omega_{s_0}=\{x|\ 0<d(x,\p\Omega)\le s_0\}$.  We have the following:

\proclaim{Theorem 4.2} Let $(\Omega^3,g)$ be a bounded domain with smooth boundary $\Sigma$ with nonnegative  scalar curvature $\cR$ such that the Gauss curvature of $\Sigma$ is positive. Let 
$$
\xi=\frac14\cR_{\min}^\frac12 
$$
where  $\cR_{\min}=\inf_{\Omega_{s_0}}\cR$. Let  $H$ be  the mean curvature of $\{x|\ d(x, \p\Omega)=s\}$, $0\le s\le s_0$, $H_+=\max\{H,0\}$ and $H_{\min}=\min_{\p\Sigma}H$. Suppose
\roster
\item"{(i)}" $\xi\ge H_+\tanh(\xi s_0)$ in $\Omega_{s_0}$; and
\item"{(ii)}" $\xi\tanh(\xi s_0)\ge -4H_{\min}$.
\endroster
Then
$$
\BY(\Omega)\ge \frac14|\Sigma|\, \xi\tanh(\xi
s_0).
$$
In particular, $\BY (\Omega)$ is bounded below by a nonnegative constant depending on $\cR_{\min}$, $H_+$, $H_{\min}$, $|\Sigma|$ and $s_0$.
\endproclaim
\proclaim{Remark 4.2} It is easy to se that if $\cR\ge0$ and $H\ge0$, then the conditions in the theorem will be satisfied. Also, the theorem says that if $\cR>0$, then the Brown-York mass of $\Omega$ is still nonnegative provided that mean curvature of its boundary is not very negative.
\endproclaim

\demo{Proof} Since $\Delta-\cR/8$ has positive first Dirichlet eigenalue, (4.2) has a unique positive solution $u$
with $q=\cR/8$. Note that $0<u\le 1$. Let $\xi$ be as in the assumptions of the theorem and let
$$
\phi(s)=\cosh(\xi s)-\tanh(\xi s_0)\sinh(\xi s),
$$
$0\le s\le s_0$. Then  on $0\le s\le s_0$, $\phi>0$ because $\tanh
\xi s\tanh (\xi s_0)<1$. $\phi(0)=1$,

$$
\frac{d^2\phi}{ds^2}=\xi^2\phi,\tag4.20
$$

$$
\frac{d\phi}{ds} =\xi\lf( \sinh(\xi s)-\tanh(\xi s_0)\cosh(\xi s)\ri) \tag4.21
$$
and so
$\frac{d\phi}{ds}<0$ in $0\le s<s_0$ and $\frac{d\phi}{ds}=0$ at $s=s_0$.

Define $f(x)=\phi\lf(d(x,\p\Omega)\ri)$ for $x\in \Omega_{s_0}$, and $f(x)=\phi(s_0)$ for $x\in \Omega\setminus \Omega_{s_0}$. Then $f$ is Lipschitz in $\Omega$ and in $\Omega_{s_0}$

$$
\split
\Delta f-\frac{\cR}8f&=\phi''-H\phi'-\frac{\cR}8 \phi\\
&\le \lf(\xi^2-\frac{\cR_{\min}}8\ri)\phi+H_+\xi\lf(\tanh(\xi s_0)\cosh(\xi s)-\sinh(\xi s)\ri)\\
&\le   \xi\cosh(\xi s)\lf[-\xi\lf(1-\tanh(\xi s_0)\tanh (\xi s)\ri)+ H_+\lf(\tanh(\xi s_0)-\tanh(\xi s)\ri)\ri]
\endsplit\tag4.22
$$
where we have used the fact that $\cR_{\min}=16\xi^2$, $\phi>0$  and $\phi'\le 0$ in $0\le s\le s_0$.

On the other hand,
$$
\split
\xi&\lf(1-\tanh(\xi s_0)\tanh(\xi s)\ri) \ge H_+\lf( \tanh(\xi s_0)-\tanh(\xi s)\ri)\\
&\iff \xi-H_+ \tanh(\xi s_0) \ge \lf( \xi \tanh(\xi
s_0)-H_+\ri)\tanh(\xi s).
\endsplit\tag4.23
$$
Since
$$
\xi+H_+\ge (\xi+H_+)\tanh(\xi s_0),
$$
which implies
$$
\xi-H_+ \tanh(\xi s_0)\ge \xi \tanh(\xi
s_0)-H_+,
$$
we have
$$
\xi-H_+\tanh(\xi s_0) \ge \lf( \xi \tanh(\xi
s_0)-H_+\ri)\tanh(\xi s) \tag4.24
$$
if $\xi\tanh\frac{\xi s_0}2-\alpha  \ge0$. The above inequality is obvious true if   $\xi\tanh\frac{\xi s_0}2-\alpha<0$
because $\xi\ge H_+\tanh(\xi s_0)$ by condition (i) in the assumptions.
From (4.22)--(4.24), we conclude that
$$
\Delta f-\frac{\cR}8f\le 0\tag2.25
$$
in $\Omega_{s_0}$. Since $\cR\ge0$, (2.25) is also true in the interior of $\Omega\setminus \Omega_{s_0}$. Since $\phi'(s_0)=0$, it is easy to see that $f$ satisfies (2.25) weakly in $\Omega$. By the maximum principle and by the fact that $f=1$ on $\p\Omega$, we conclude that $u\le f$ in $\Omega$. Hence

$$
\frac{\p u}{\p\nu}\ge \frac{\p f}{\p\nu}=\xi  \tanh(\xi s_0)
\tag4.26
$$
on $\p\Omega$, where $\nu$ is  the unit outward normal of $\p\Omega$.

Consider the metric $g_1=u^4g$, then the scalar curvature of $g_1$
is zero and the mean curvature $\ol H$ with respect to $g_1$ of
$\p\Omega$ satisfies

$$
\ol H=H+\frac14\frac{\p u}{\p\nu}\ge H_{\min}+\frac14\xi\tanh(\xi
s_0)\ge0 
$$
by condition (ii) in the assumptions of the theorem. Since $u=1$ on the boundary, the induced metric on $\Sigma$ is the same as before and so the Gauss curvature of $\Sigma$ is positive. By Proposition 4.1, we have
$$
\int_{\Sigma}\lf(H_0-\ol H\ri)d\sigma\ge0,
$$
and hence
$$
 \int_{\Sigma}\lf(H_0-H\ri)d\sigma\ge|\Sigma|\cdot\frac14\xi\tanh(\xi
s_0).
$$
From this the theorem follows.

\enddemo

Similar to Corollary 4.2, by Theorem 4.2 and Lemma 4.2, we have:

\proclaim{Corollary 4.3}   With the assumptions and same notations as in Theorem 4.2. Then $\BY(\Omega)$ is bounded below as in (4.18) and (4.19) with $\cR_-=0$.
\endproclaim
\demo{Proof} By the proof of the theorem, with the same notations as in the proof of the theorem, we have
$$
\BY(\Omega)\ge\frac14\int_{\Sigma}\frac{\p u}{\p \nu}.
$$
The corollary follows as in the proof of Corollary 4.2.
\enddemo
\subheading{\S5 Some examples and applications}\vskip .2cm

In this section, we will give some examples which are related to results in previous sections. Some of the examples
might be well-known.\vskip .2cm

\noindent{\bf Example 1}: In Corollary 3.1, it is proved that if
the negative part of the scalar curvature is small compared with
the Sobolev constant and the positive part of the scalar
curvature, then the ADM mass of a spin AF manifold is nonnegative.
The following example show that if we only assume that the
negative part of scalar curvature is small compared with the
Sobolev constant, the ADM mass might still be negative.

 Let $g_{ij}=u^4\delta_{ij}$ be a conformal metric on $\R^3$. Then the scalar curvature of $g$ is

$$
\cR=-8u^{-5}\Delta_0 u, \tag5.1
$$
where $\Delta_0$ is the Euclidean Laplacian.

Let $v$ be a nonconstant smooth function such that $\Delta_0v\ge0$ in $B(1)=\{x\in \R^3|\ |x|<1\}$, say,
and $\Delta_0v=0$ outside $B(1)$ such that $v\to0$ near infinity. Then $v\le0$. We may also assume that $v>-1$.
Then $v=-\frac A{|x|}+O(|x|^{-2})$ near infinity with $A>0$. For any $1>\e>0$, let $u=1+\e v$ and consider
the metric $g_{ij}=u^4\delta_{ij}$. Then the scalar curvature $\cR\le0$ by (5.1). Moreover, if   $M=(\R^3,g)$,
then

$$
\lf(\int_M \cR_-^\frac 32 dV_g\ri)^\frac23=8\e\lf(\int_M
\lf(u^{-5}\Delta_0v\ri)^\frac32u^6dV_e\ri)^\frac23.
$$
where $dV_e$ is the Euclidean volume form. If $\e$ is small
enough, then $u$ is close to 1. Hence $\lf(\int_M R_-^\frac 32
dV_g\ri)^\frac23$ can be made arbitrarily small compared with the
Sobolev constant $\Lambda$ of $M$ by letting $\e\to0$.  But the
mass of $M$ is negative.\vskip .2cm

\noindent{\bf Example 2}: It is easy to see that assumptions of Corollary 3.1 or the assumptions in
\cite{ZZ, Theorem 4.1} imply $\int_M\cR\ge0$ for an AF manifold $M$. However, it is not hard to construct
examples of AF metrics $g$ on $\R^3$ such that $\int_M \cR>0$ but the ADM mass $\frak m_g$ of $M$ is negative,
where $M=(\R^3,g)$.

Let $g_{ij}=u^4\delta_{ij}$ be an AF metric on $\R^3$ with $\frak m_g<0$, for example the metric in the example
in \cite{ZZ}. Let $v\ge0$ be a smooth function with support in $B(1)=\{|x|<1\}$ and $v\not\equiv0$. Define
$\varphi  =1+a v$ where $a>0$. Consider the metric $\wt  g_{ij}=\lf(\varphi u\ri)^4\delta_{ij}$. Then
 $\varphi u=u$ outside $B(1)$. Hence $\frak m_{\wt g}=\frak m_g<0$. The scalar curvature $\cR_{\wt g}$ is given by

$$
\cR_{\wt g}=-\lf(\varphi u\ri)^{-5}\Delta_0(\varphi u).
$$

Hence

$$
\split
\int_{\R^3} \cR_{\wt g}dV_{\wt g}&=\int_{|x|>1}\cR_g dV_g-\int_{|x|<1}\varphi u\Delta_0\lf(\varphi u\ri)dV_e\\
&=\int_{|x|>1}\cR_g dV_g+\int_{|x|<1}|\nabla_0 (\varphi u)|^2dV_e-\int_{|x|=1}u\frac{\p u}{\p r}\\
&= \int_{|x|>1}\cR_g dV_g-\int_{|x|=1}u\frac{\p u}{\p r}+\int_{|x|<1}|\nabla_0u+a\nabla (vu))|^2dV_e\\
&\ge\int_{|x|>1}\cR_g dV_g-\int_{|x|=1}u\frac{\p u}{\p
r}-\int_{|x|<1}|\nabla_0u|^2+\frac12 a^2\int_{|x|<1}
|\nabla_0(vu)|^2dV_e.\endsplit
$$
Choose $v$ so that $vu$ is not constant in $|x|<1$ and choose $a$ large enough, we have $\int_{\R^3}\cR_{\wt g}dV_{\wt g}>0$.

\vskip .2cm

In  \cite{B1, Theorem 5.2}, it was proved that if a metric $g$ is close to the Euclidean metric of $\R^3$ and if $\int_{\R^3}\cR_g dV_e\ge0$, then $\frak m_g\ge0$. On the other hand, we have the following observation.

\proclaim{Proposition 5.1} Suppose $g_{ij}=u^4\delta_{ij}$ is an AF metric on $\R^3$. Then $\frak m_g\le C\int_M\cR dV_g$ for some   absolute constant $C>0$. Equality holds if and only if $g$ is Euclidean. In particular, if $\frak m_g\ge0$, then $\int_{\R^3} \cR dV_g\ge0$.

\endproclaim

\demo{Proof}   Let $v=1-u$. Then $v=\frac Ar+O(r^{-2})$, and $|\partial v|=O(r^{-2})$ etc. Then the mass of $g$ is given by

$$
C\frak m_g=- \int_{\p B(\infty)}\frac{\p u}{\p
r}=-\int_{\R^3}\Delta_0 udV_e \tag5.2$$
for some absolute constant $C>0$.
Note that

 $$
\split
  \int_M \frac \cR8dV_g
&=-\int_{\R^3} \lf(u^{-5}\Delta_0 u\ri)u^6dV_e\\
&=-\int_{\R^3}u\Delta_0 udV_e
\endsplit\tag5.3
$$
By (5.2) and (5.3), we have
$$
\split C\frak m_g-\int_M \frac \cR8dV_g
&=\int_{\R^3}\lf(-\Delta_0u+u\Delta_0u\ri)dV_e\\
&=\int_{\R^3}(-1+u)\Delta_0udV_e\\
&=\int_{\R^3}v\Delta_0 vdV_e\\
&=-\int_{\R^3}|\nabla_0 v|^2dV_s\\
&\le 0
\endsplit
$$
because $v=O(r^{-1})$ and $|\nabla_0 v|=O(r^{-2})$. From this it is easy to see the proposition is true.\enddemo

\noindent{\bf Example 3}: There are examples of AF metrics defined on $\R^3$ with zero scalar curvature and positive ADM mass. In fact, Miao \cite{M} constructed an AF metric on $\R^3$ which is scalar flat, conformally flat outside a compact set and contains a horizon. In particular, the mass is positive.

We may also construct the scalar flat but nonflat AF metrics in
the following way.  Take a metric $g$ on $\Bbb S^3$ so that $(\Bbb
S^3, g)$   is not conformal to the standard metric and has
positive Yamabe invariant.  Take a point p in $\Bbb S^3$ and
consider the metric $u^4 g$ on $\Bbb S^3\setminus\{p\}$, where $u$
is the Green's function for the conformal Laplacian with pole at
$p$, which exists and positive by \cite{LP, \S6}. Then the
manifold $M=\lf(\Bbb  S^3\setminus\{p\}, u^4 g\ri)$ is AF, scalar
flat, and $M$ is diffeomorphic to $\R^3$. Note that $M$ has
positive mass and M is not conformal to $\R^3$, see \cite{LP,
\S11}. To construct $g$, we can perturb the standard metric in
some neighborhood of a point $p$ so that it is not conformally
flat in that neighborhood.   One may assume the perturbation is
small so that the scalar curvature is still positive. Then the
Yamabe invariant of the metric must be nonnegative by the
definition of the Yamabe functional. It must be positive,
otherwise we can find a positive solution of conformal Laplacian.
The solution must be   constant by the strong maximum principle,
which is impossible because the scalar curvature of $g$ is
positive. The metric is not conformal to the standard metric
because it is not locally conformally flat. So we are done.\vskip
.2cm

\noindent{\bf Example 4}: From Example 3, we can construct AF metrics $g$ on $\R^3$ so that $\cR\le 0$,
$\cR\not\equiv0$ and $\frak m_g>0$.

To do this, let $g$ be an AF metric on $\R^3$   defined in Example 3 so that $\cR_g\equiv0$ and $\frak m_g>0$.
Let $v$ be a nonconstant bounded subharmonic function with respect to $g$ such that $v$ is harmonic outside a
compact set so that $v\to0$ near infinity. Then $v\sim a/r$ near infinity for some constant $a$. Let $u=1+\e v$
where $\e>0$ is small and let $\wt g=u^4 g$. Then $\wt g$ is AF metric on $\R^3$ and there is an absolute constant
$C>0$ such that

$$
\frak m_{\wt g}=-C\int_{\Bbb S_\infty}\frac{\p u}{\p r}+\frak m_g=
-C\e \int_{\Bbb S_\infty}\frac{\p v}{\p r}+\frak m_g.
$$
Hence $\frak m_{\wt g}>0$ if $\e $ is small enough. On the other hand, since $\cR_g\equiv0$,
$$
\cR_{\wt g}=-u^{-5}\Delta_g u=-\e u^{-5}\Delta_g v\le0
$$
and $\cR_{\wt g}<0$ somewhere because $\Delta_g v>0$ somewhere. In particular,
$\int_{\R^3}\cR_{\wt g}dV_{\wt g}<0$.\vskip .2cm

\noindent{\bf Applications}: We now discuss some relations of the results in \cite{ST1}  and
\cite{H} and the classical Minkowski's inequalities for convex bodies in $\R^3$.

It was proved  by Herzlich  \cite{H, Proposition 2.1}, that if $(M,g) $ is an AF 3-dimensional manifold with nonnegative scalar curvature with an inner boundary $\Sigma$ which is homeomorphic to $\Bbb S^2$ whose mean curvature with respect to the inner normal satisfies

$$
H\le 4\sqrt{\frac{\pi}{A(\Sigma)}}\tag5.4
$$
where $A(\Sigma)$ is the area of $\Sigma$, then the mass of $M$ is nonnegative. By convention the mean curvature of $\R^3\setminus B(1)$ in the Euclidean space is 2. On the other hand, the well-known Minkowski's inequalities for convex bodies in $\R^3$ state that if $\Sigma$ is a compact   convex surface in $\R^3$, then
$$
\lf(\int_\Sigma H_0\ri)^2\ge 16\pi A,\tag5.5
$$
and
$$
\frac{4A^4}{9V^2}\ge \lf(\int_\Sigma H_0\ri)^2 \tag5.6
$$
where $H_0$ is the mean curvature of $\Sigma$ in $\R^3$, $A$ is the area of $\Sigma$ and $V$ is the volume of the region bounded by $\Sigma $, see \cite{BG, p. 438}. Moreover equality holds either in (5.5) or (5.6) if and only if $\Sigma$ is a standard sphere.

 Using the results in \cite{H,ST1}, one can derive (5.5).  In fact, using the method of \cite{ST1}, one can find an AF metric $g$ on the exterior $\R^3_\Sigma$ of $\Sigma$ in $\R^3$ with zero scalar curvature such that   $\Sigma$ has constant mean curvature  $H\equiv 4\sqrt{\frac{\pi}{A(\Sigma)}}$. Note that by the construction in \cite{ST1}, the metric $g$ when restricted on $\Sigma$ is the same as the restriction of the Euclidean metric. By \cite{H}, $\frak m_g\ge0$. Hence by the result of \cite{ST1},

$$
\int_\Sigma(H_0-H)\ge0.
$$
This implies (5.5).

We can also use (5.6) and the result of \cite{ST1} to prove that
condition (5.4) is sharp in  the result of \cite{H} mentioned
above in the following sense. Given any $\e>0$, we can find a
manifold with boundary satisfying all the conditions of
Proposition 2.1 in \cite{H}  except that
$$
H\ge 4\sqrt{\frac\pi{A(\Sigma)}}+\e
$$
and the mass of the manifold is negative. To construct such an example, note that equalities hold  in (5.5) and (5.6) if and only if $\Sigma$ is the standard sphere. Hence for any $\e>0$, we can perturb the standard sphere so that it is still strictly convex, but
$$
\frac{4A^4}{9V^2}< 16\pi A+\e.
$$
Now we find the metric $g$ in the exterior of $\Sigma$ as before with initial mean curvature $H$ such that

$$
\lf(\int_\Sigma H\ri)^2=16\pi A+\e.
$$
If the mass is nonnegative, then we have
$$
\lf(\int_\Sigma H_0\ri)^2\ge \lf(\int_\Sigma H\ri)^2=16\pi
A+\e>\frac{4A^4}{9V^2}
$$
by \cite{ST1}, which is impossible because of (5.6).

\enddocument
\end